\documentclass[reqno,12pt]{amsart}

\usepackage[all]{xy}
\usepackage{graphicx,enumerate}

\numberwithin{equation}{section}

\theoremstyle{plain}
\newtheorem{lemma}{Lemma}[section]
\newtheorem{theor}[lemma]{Theorem}
\newtheorem{prop}[lemma]{Proposition}

\newtheorem{claim}[lemma]{Claim}
\newtheorem{claim*}{Claim}

\newenvironment{concl}{\begin{bf}Conclusion.\end{bf}\thinspace  \begin{em}}{\end{em}}

\theoremstyle{definition}

\theoremstyle{remark}
\newtheorem{rem}[lemma]{Remark}

\newcommand{\pic}{\operatorname{Pic}} 
\newcommand{\Pic}{\operatorname{Pic}}
\newcommand{\bs}{\operatorname{Bs}}

\newcommand{\edim}{\operatorname{edim}}

\newcommand{\p}{\mathbb{P}}

\newcommand{\z}{\mathbb{Z}}

\newcommand{\co}{{\mathcal O}}
\newcommand{\cl}{{\mathcal L}}

\newcommand{\ca}{{\mathcal A}}

\newcommand{\cs}{{\mathcal S}}

\newcommand{\ci}{{\mathcal I}}
\newcommand{\cw}{{\mathcal W}}

\newcommand{\cx}{{\mathcal X}}

\newcommand{\nvs}{\vspace{-2mm}}

\newcommand{\map}{\rightarrow}
\newcommand{\rmap}{\dashrightarrow}

\begin{document}

\title{Linear systems on a class of anticanonical rational threefolds}
\author{Cindy De Volder}
\address{
Department of Pure Mathematics and Computeralgebra,
Krijgslaan 281, S22, \newline B-9000 Ghent, Belgium}
\email{cindy.devolder@ugent.be}

\author{Antonio Laface}
\address{Departamento de Matematicas, Universidad de Concepci\'on. \newline
Casilla 160-C Concepci\'on, Chile}
\email{antonio.laface@gmail.com}

\begin{abstract}
Let $X$ be the blow-up of the three dimensional complex projective space along $r$ general points of a smooth elliptic quartic curve $B\subset\p^3$ and let $L\in\Pic(X)$ be any line bundle. The aim of this paper is to provide an explicit algorithm for determining the dimension of $H^0(X,L)$.
\end{abstract}
\maketitle


\section{Introduction}
Let $X$ be the blow-up of the three dimensional projective space along $r$ general points lying on a smooth elliptic quartic curve $B$.
The aim of this paper is to provide an explicit algorithm for determining the dimension of
$H^0(X,L)$
for any $L\in\pic(X)$. 
This dimension depends of course on the degree and multiplicities of the general divisor of $\p^3$ corresponding to $L$. In this paper we show that in fact this number is completely determined by the values of the intersections $l_i\cdot L$ and $C\cdot L$, where the $l_i$'s are the strict transforms of lines through pair of the $r$ points and $C$ is the strict transform of $B$.

This work is an attempt to generalize the results of~\cite{bh1} to the three dimensional case by extending the techniques used in~\cite{DL8}.
Recently, in~\cite{CT} and~\cite{CEPG}, an higher dimensional analog of the same problem is studied under the more restrictive hypothesis that all the points lie on a rational normal curve of $\p^n$. 
It turns out that this assumption implies the finite generation of the Cox ring of the blow-up variety, while in the case analyzed by the present paper this statement is false.

The paper is organized as follows: in section~\ref{notation} we fix the necessary notation while section~\ref{prel}
focuses on preliminary results regarding the intersection theory of the varieties which are needed throughout the paper.
The main algorithm is explained in section~\ref{algorithm}, here we  show how starting from the linear system $\cl$, associated to $L\in\pic(X)$, one can find a fixed-component free linear system $\cl'$ of the same dimension of $\cl$. Then we proceed to define three different types of systems, listed in conclusion~\ref{concl 3 types}, which cover the range of all the possibilities.
The dimension of a linear system in each one of these classes is given explicitly in theorems~\ref{t=0}, \ref{2m,m^r} and \ref{t>0,a}.

\newpage

\section{Notation}\label{notation}

The aim of this section is to provide the necessary notation for linear systems defined on blow-ups of $\p^2$ and $\p^3$ and a quadric.
In what follows the ground field is assumed to be algebraically closed of characteristic 0. \\

\begin{center}
\begin{figure}[htp]
\includegraphics[scale=.5]{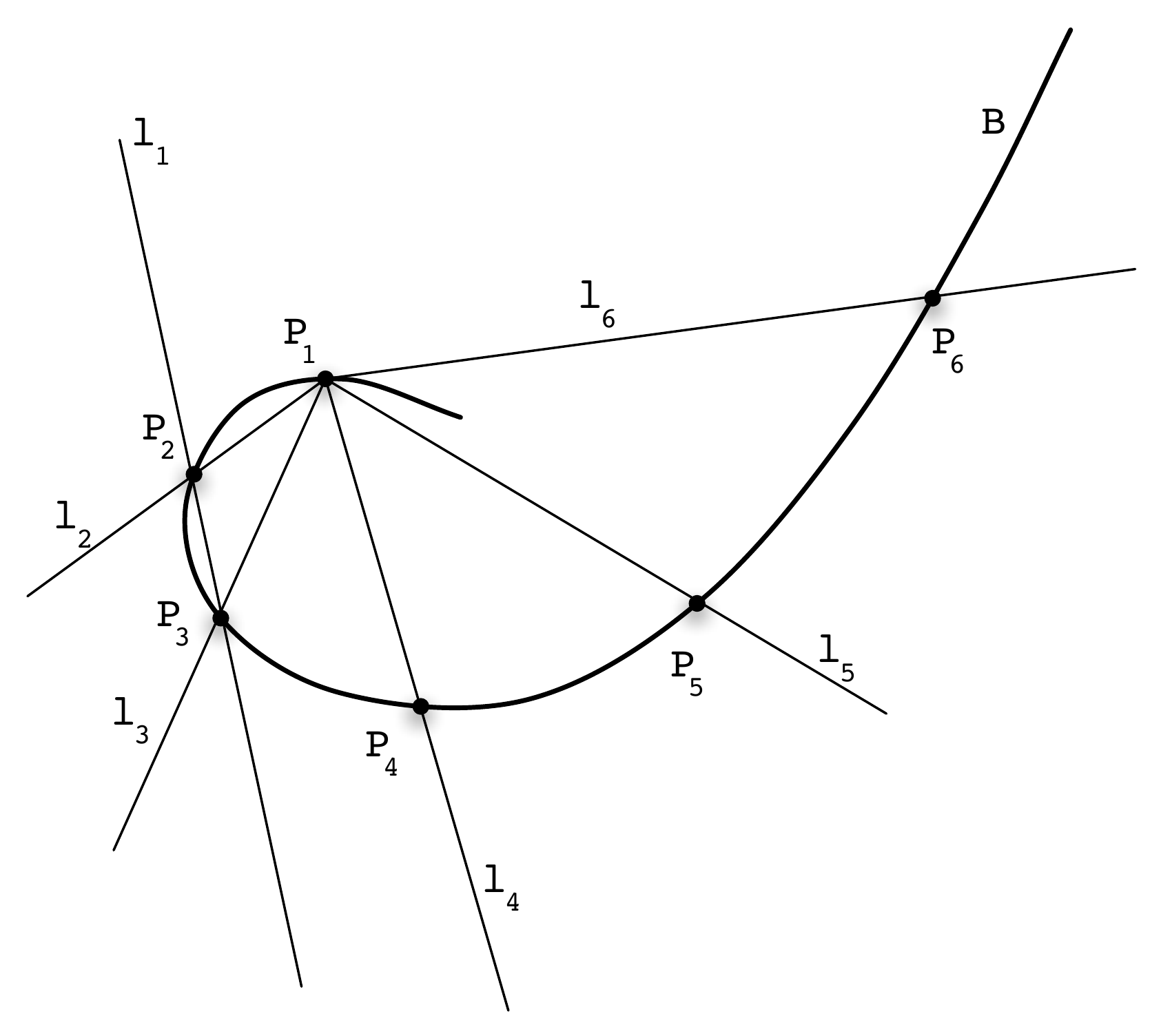}
\caption{The configuration of points and lines}
\end{figure}
\end{center}

{\bf The definition of $X$}
We start with a smooth quadric $Q\subset\p^3$ and a general $B\in |-K_Q|$, i.e. $B$ is an elliptic curve of degree $4$.
On this curve we choose $p_1, \ldots, p_r$ points in general position
and $Z = m_1 p_1 + \dots + m_r p_r$ is a zero-dimensional subscheme with $m_i$ non-negative integers and associated ideal sheaf $\ci_Z$. With abuse of notation we denote by  
\[
\cl_3(d;m_1,\ldots,m_r)
\]
both the sheaf $\co_{\p^3}(d) \otimes \ci_Z$ and its associated linear system.

The {\em expected dimension} of such a linear system $\cl$ is $\edim(\cl) := \max\{ -1, v(\cl) \}$, where
\[
v(\cl) = \binom{d+3}{3} - \sum_{i=1}^r \binom{m_i+2}{3} -1,
\]
is the {\em virtual dimension}.

Let $\pi : X \map \p^3$ be the blow-up map of $\p^3$ along $p_1,\ldots, p_r$, then we will denote by $H$ the pull-back of a plane and by $E_i$ the exceptional divisor corresponding to $p_i$.
In this way the linear system $\cl_3(d;m_1,\ldots,m_r)$ will correspond to the complete linear system
\[
\cl_{X}(d;m_1,\ldots,m_r) := |dH - m_1 E_1 - \cdots - m_r E_r|
\] 
defined on $X$. As before we will use the same notation for the linear system and its associated invertible sheaf.
We will say that a system, on $\p^3$ or on $X$, is {\em standard} if
\[
m_1 \geq \ldots \geq m_r \geq 0\hspace{1cm} {\rm and}\hspace{1cm} 2d \geq m_1 + \cdots + m_4,
\]
while we we call it {\em almost standard} if it becomes standard after reordering its multiplicities. The strict transform of $B$ will be denoted by $C$. \\

{\bf The blow-up of $Q$} We have the following commutative diagram
\[
\xymatrix{
Q_r\ar[r]\ar[d]_{\pi_{Q}} & X\ar[d]^{\pi} \\
Q\ar[r] & \p^3
}
\]
where $\pi_Q$ is the blow-up map of $Q$ at the $p_i$'s and
$e_i := E_i\cap Q_r$
is the exceptional divisor on $Q_r$ corresponding to $p_i$.
We will use the notation 
\[
\cl_Q(a,b;m_1,\ldots,m_r) := |ah_1+bh_2-m_1 e_1 - \ldots - m_r e_r|
\]
where $ah_1+bh_2 = \pi_Q^*~ \co(a,b)$
to denote both the linear system and its corresponding invertible sheaf. 
We will say that such a linear system is in {\em standard form}
if 
\[
m_1 \geq \ldots \geq m_r \geq 0\hspace{1cm} {\rm and}\hspace{1cm} a+b \geq m_1 + m_2 + m_3 + m_4
\]
and it is {\em standard} if there exists a base of 
$\Pic(Q_r)$ such that the system is standard when written in that base.

A proof of the following proposition,which will be used several times in this paper, can be easily obtained by readapting the arguments of~\cite{bh1} to the blow-up of a quadric along points on a smooth anticanonical divisor.

\begin{prop}\label{har}
Let $D$ be an effective divisor of $Q_r$, and let $C\in |-K_{Q_r}|$ be a smooth curve, then $C\subseteq\bs |D|$ if and only if $D\cdot C \leq 0$.
\end{prop}

{\bf The definition of $Y_I$}
Consider the strict transform $l_i$ of the line $\overline{p_1p_i}$ and denote by $l_1$ the strict transform of the line through $p_2$ and $p_3$.
Given a subset $I \subset \{1,\ldots,r \}$, define 
\[
\pi_I : Y_{rI} \map X
\]
to be the blow-up of $X$ along all the $l_i, i\in I$ and denote by
$F_i$, the exceptional divisor corresponding to $l_i$.

When there is no ambiguity about $r$, we will use $Y_I$ instead of $Y_{rI}$.

By abuse of notation, we will use $H$, resp. $E_i$, on $Y_{I}$, 
to denote the pull-back of $H$, resp. $E_i$;
and we let $C$, resp. $\tilde{Q}$ to denote the corresponding strict transforms.
With
\[
\cl_{Y}(d; m_1,\ldots,m_r; \{ t_i \}_{i \in I}) := |dH - m_1E_1 - \ldots -m_r E_r - \sum_{i\in I} t_i F_i|,
\]
where with $d , m_i , t_i \geq 0$, we will denote the complete linear system and its corresponding invertible sheaf. \\

{\bf The definition of $\tilde{Y}_I$}
Our last object will be the blow-up of $Y_I$ along $C$ defined as 
\[
\tilde{\pi}_I : \tilde{Y}_I \map Y_{I}.
\]
The exceptional divisor over $C$ will be denoted by $F$ and
as before, by abuse of notation, we will denote the pull-back of resp. $H$, $E_i$ and $F_j$ with the same letters and the same for the strict transform of $\tilde{Q}$.
With
\[
\cl_{\tilde{Y}}(d; m_1,\ldots,m_r; \{ t_i \}_{i \in I}; t) := |dH - m_1E_1 - \ldots -m_r E_r - \sum_{i\in I} t_i F_i - tF|,
\]
where $d , m_i , t_i, t \geq 0$, we will denote the complete linear system and its corresponding invertible sheaf.

In all the above notation, if $m_i = m_{i+1} = \cdots = m_{i+k} = m$, then we will
write $m^k$ in stead of $m_{i},m_{i+1},\ldots,m_{i+k}$.

\section{Preliminaries}\label{prel}

In this section we deal with the varieties $X, Y_I$ and $\tilde{Y}_I$ just defined and we will work out their Chern classes and the intersection product of cycles on them.

\vspace{5mm}\subsection{A Cremona transformation of type (3,3)}~

To any four non-collinear points of $\p^3$ we can associate a $(3,3)$ birational map corresponding to the linear system $\cl_3(3;2^4)$. After a linear change of coordinates this map can be described by:
\begin{equation}\label{cubic}
\phi(x_0:x_1:x_2:x_3) =(x_0^{-1}:x_1^{-1}:x_2^{-1}:x_3^{-1}).
\end{equation}
This birational map induces a quasi-isomorphism $\tilde{\phi}: X\rmap X'$ of the blowing-up of two $\p^3$'s along four points.
As a consequence, the induced map $\phi^*: \Pic(X')\rightarrow\Pic(X)$ is an isomorphism which is described by the following:
\begin{prop}[\cite{lu}]\label{cre-surfaces}
The action of transformation~(\ref{cubic}) on
$\cl=\allowbreak\cl_X(d;\allowbreak m_1,\allowbreak \ldots,m_4)$
is given by:
\begin{eqnarray}\label{cre-a1}
\phi^*(\cl) & := & \cl_X(d+k;m_1+k,\ldots, m_4+k),
\end{eqnarray}
where $k=2d-\sum_{i=1}^4m_i$.
\end{prop}
Note that the points $\phi(p_1),\ldots,\phi(p_4)$ are still in general position $\phi(C)$, therefore the action of $\phi$ on a linear system $\cl_{3}(d;m_1,\ldots,m_r)$ results in a linear system
$\cl_{3}(d+k;m_1+k,\ldots, \allowbreak m_4+k,\allowbreak  m_5, \ldots,m_r)$ and the same is true for the corresponding systems on $X$.

Observe that even if $\dim \phi^*(\cl) = \dim \cl$, in general the
virtual dimensions of the two systems may be different (see e.g.~\cite{lu}). 

\vspace{5mm}\subsection{Euler characteristic and Chern classes of the fundamental varieties}~

Let $D$ be a divisor on a smooth rational threefold. According to the Riemann-Roch formula, the Euler characteristic of $\co_X(D)$ is
\[
\cx (\co_X(D)) = \frac{1}{12}[D (D - K_X) (2D - K_X) + c_2(X) D] + 1,
\]
where $c_2(X)$ is the second Chern class of $X$ and $K_X$ is the canonical class of $X$.

The threefolds we are going to work on will be $\p^3$, $X$, $Y_{I}$ and $\tilde{Y}_{I}$. Hence we need to know $K_X$, $c_2(X)$ and for these varieties.

\begin{prop}[see e.g.~\cite{GH}]\label{c's p1}
Let $\pi: \tilde{X} \map X$ be the blow-up of a smooth algebraic threefold $X$ along a smooth, connected subscheme $V\subset X$ of dimension $a\leq 1$ and let $E$ be the exceptional divisor, then:
\begin{align*}
c_1(\tilde{X}) &= \pi^* c_1(X) + (a-2)~E \mbox{ and} \\
c_2(\tilde{X}) &= \pi^* (c_2(X) + a~ \eta_V) - a~ \pi^* c_1(X)~E,
\end{align*}
where $\eta_V\in H^4(X,\z)$ is the class of the curve $V$ in $X$.
\end{prop}

Using the fact that $c_1(\p^3) = 4H$ and $c_2(\p^3) = 6 H^2$, we can deduce the following:

\begin{lemma}\label{c1,c2}
Let $I\subseteq\{1,\ldots,r\}$ be a subset of cardinality $a$ and 
denote by $\epsilon$ a number which is equal to $1$
if $1\in I$ and $0$ otherwise.
The first two Chern classes of $X, Y_I$ and $\tilde{Y}_I$ are:
\begin{itemize}\itemsep 5mm
\item[$X$] 
	\begin{itemize}\itemsep 2mm
		\item[] $c_1 = 4H - 2E_1 - \cdots - 2E_r$
		\item[] $c_2 = 6H^2$.
	\end{itemize}
\item[$Y_I$]
	\begin{itemize}\itemsep 2mm
		\item[] $c_1 = 4H - 2E_1 - \cdots - 2E_r - \sum_{i \in I} F_i$ 
		\item[] $c_2 = (6+a)H^2 + aE_1^2 + \sum_{i\in I\setminus\{1\}} E_i^2 + \epsilon~ (E_2^2 + E_3^2 - E_1^2)$
	\end{itemize}
\item[$\tilde{Y}_I$]
	\begin{itemize}\itemsep 2mm
		\item[] $c_1 = 4H - 2E_1 - \cdots - 2E_r - \sum_{i \in I} F_i - F$
		\item[] $c_2 = \pi^* c_2(Y_I) + 4H^2 + E_1^2 + \cdots + E_r^2  - 4HF + 2E_1F + \cdots + \nvs 2 E_rF$.
	\end{itemize}
\end{itemize}
\end{lemma}

\begin{proof}
The first Chern class of all these varieties is easily determined by means of Proposition~\ref{c's p1} and the same is true for 
$c_2(X)$.

To obtain $c_2(Y_I)$, we will assume for simplicity that $I = \{2,\ldots,s\}$.
In this case, we can consider $\pi_I$ as the composition
$\pi_s \circ \pi_{s-1} \circ \cdots \circ \pi_2$, where 
$\pi_i : Y_{I_i} \map Y_{I_{i-1}}$ the blow-up map of $l_i$ and $Y_{I_1} = X$.
Since the class of $l_2$ in the Chow ring of $X$ is given by
\[
\eta_{l_2} = (H - E_1 - E_2)^2 = H^2 + E_1^2 + E_2^2,
\]
then, by applying Proposition~\ref{c's p1}, we obtain the following:
$$
c_2(Y_{I_2}) = 6H^2 + H^2 + E_1^2 + E_2^2 - (4H - 2E_1 - \cdots - 2E_r)~F_2.
$$
Now observe that $(4H - 2E_1 -\dots -2E_r)~l_2 = 0$ in $X$ so that on $Y_I$ we have:
\[
(4H - 2E_1 -\dots -2E_r)~F_2 \equiv 0.
\]
Repeating this argument for all $i = 2,\ldots, s$,
and using the fact that $F_i ~F_j \equiv 0$ for all $i \neq j$ (because the curves $l_i$ and $l_j$ do not intersect),
we obtain that
\begin{align*}
c_2(Y_{I_2}) &= 6H^2 + \sum_{i=2}^s (H^2 + E_1^2 + E_i^2)
\end{align*}
which gives the desired result.

For determining $c_2(\tilde{Y}_I)$ we need to know
the class of $C$ in $Y_I$, which is given by:
\[
\eta_{C} = (2H - E_1 - \cdots -E_r)^2 = 4H^2 + E_1^2 + \cdots + E_r^2.
\]
From Proposition~\ref{c's p1} and $F_i~ F \equiv 0$ (because $l_i$ and $C$ do not intersect) we deduce:
\begin{align*}
c_2(\tilde{Y}_I) & = \pi^* c_2(Y_I) + 4H^2 + E_1^2 + 
	\cdots + E_r^2 - \pi^*c_1(Y_I)~ F 
		\end{align*}
which gives the desired result.
\end{proof}

\vspace{5mm}\subsection{Intersection products}~

The aim of this section is to provide explicit formulas for the intersection of divisors on $X$, $Y_I$ and $\tilde{Y}_I$. In particular, given any three divisors on one of these varieties, we are interested in evaluating their intersection number.
Since this number and that of the intersection of their pull-backs are the same, we need only to work out the calculations for divisors in $\tilde{Y}$.

The following proposition will be needed in what follows.

\begin{prop}[see e.g. \cite{GH}]\label{E|_E}
Let $X$ be a smooth irreducible threefold and $S \subset X$ a smooth irreducible surface, then $c_1(X)|_S = c_1(S) + S|_S$.
\end{prop}

\begin{lemma}\label{int Y_I}
The non-vanishing intersection products of any three divisors (or their pull-backs) on $X, Y_I$ and $\tilde{Y}_I$ are given below:
\begin{itemize}\itemsep 5mm\leftskip 10mm
\item[$X$\hspace{1cm}] 	$H^3 = E_i^3 = 1$ 
\item[$Y_I$\hspace{1cm}] $(H\ |\ E_2\ |\ E_3)~ F_1^2 = (H\ |\ E_1\ |\ E_i)~ F_i^2 = -1\hspace{5mm} ( i=2,\dots,r )$ \\

\noindent
$F_i^3 = 2\hspace{5mm} ( i=1,\dots,r )$
\item[$\tilde{Y}_I$\hspace{1cm} ]$H F^2 = -4$\\ 

\noindent
$E_i F^2 = -1\hspace{5mm} ( i=1,\dots,r )$\\ 

\noindent
$F^3 = 2(r-8)$
\end{itemize}

\end{lemma}

\begin{proof}
We begin by observing that the following intersections ( $i\neq j$ ):
\[
H~ E_i,\ E_i~ E_j,\ F_i~ F_j,\ F_i~ F,\ E_1~ F_1 
\]
are numerically equivalent to $0$ and the same is true for
$E_i~ F_1$ if $i \neq 2,3$ and $E_i~ F_j$ if $i,j \neq 1$ and $i \neq j$.
So, any intersection product of three divisors containing one of the above monomials must vanish.

Since $H^2$ is the pull-back of the class of a line of $\p^3$, it has non-zero intersection only with $H$ and this intersection is clearly equal to $1$. 
A general $W \in |H|$ is the blow-up of a plane along $\# I + 4$ points (because a general $H$ of $X$ intersects every $l_i$ in one point and $C$ in 4 points).
On $W$, let $h$ be the class of a line, $e_i$ the exceptional curve coming from the intersection point with $l_i$, and $c_1, c_2, c_3, c_4$ those coming from the intersection points with $C$. Then
$H F_i^2 = e_i^2 = -1$ and $H F^2 = (c_1  + c_2 +c_3 +c_4)^2 = -4$.

Observe that, on $X$, one has that $-E_i^2$ is the class of a line of $E_i$, hence the same is true when we consider the pull-back of this class to $\tilde{Y}_I$. This implies that $E_i^2$ has non-zero intersection only with $E_i$ and this intersection is $1$. 

If $E_i$ and $F_j$ intersect, it is easy to see that this intersection is a $(-1)$-curve on $F_j$, hence we obtain that
$E_i F_j^2 = -1$ and in the same way we have that $E_i F^2 = -1$.

On $F_i$, which is a ruled rational surface,
let $h$ denote the class of a section and $f$ the class of a fiber.
Then $c_1(F_i) \equiv 2h + 2f$,
and it follows from Proposition~\ref{E|_E} that
\begin{align*}
F_i|_{F_i} & \equiv (4H - 2E_1 - \cdots - 2E_r - \sum_{j \in I} F_j - F)|_{F_i} -  2h - 2f \\
	& \equiv 4f - 2f - 2f - {F_i}|_{F_i} - 2h - 2f \\
	& \equiv - {F_i}|_{F_i} - 2h - 2f.
\end{align*}
So, $F_i|_{F_i} \equiv -h -f$, which implies that
$F_i^3 = (-h-f)^2 = 2$.

Similarly, on $F$, which is an elliptic ruled surface,
let $h$ denote the class of a section and $f$ the (numerical equivalence) class of a fiber.
Then $c_1(F) \equiv 2h$, and from Proposition~\ref{E|_E} we obtain that
\begin{align*}
F|_{F} & \equiv (4H - 2E_1 - \cdots - 2E_r - \sum_{j \in I} F_j - F)|_{F} -  2h \\
	& \equiv 16f - \sum_{i=1}^r(2f) - F|_{F} - 2h \\
	& \equiv - F|_{F} - 2h - 2(r-8)f.
\end{align*}
So, $F|_{F} \equiv -h -(r-8)f$, which implies that
$F^3 = ( -h -(r-8)f)^2 = 2(r-8)$.
\end{proof}

\vspace{5mm}\subsection{Cohomology of an invertible sheaf and its pull-back}~

The aim of this section is to recall a criterion for comparing the cohomology of line bundles on a smooth projective variety $M$ together with the cohomology of line bundles on its blow-up $\tilde{M}$ along a smooth subscheme $V\subset M$. To this purpose consider the blow-up map $\pi: \tilde{M}\rightarrow M$ of $V\subset M$, then we have the following:

\begin{lemma}\label{H^i inv sheaf}
For any positive integer $i$ we have that $R^i \pi_* \co_{\tilde{M}} = 0$, which in turn implies:
$$H^i(\tilde{M},\pi^* \cl) \cong H^i(M,\cl)$$
for any $\cl\in\pic(M)$.
\end{lemma}

\begin{proof}
The vanishing of the higher direct images of $\pi_* \co_{\tilde{M}}$ 
depends on the fact that, for $p\in M$, the fiber
$F_p := \pi^{-1}(p)$ is either a point or $\p^r$, where $r$ is the codimension of $V$ in $M$.  
From the projection formula~\cite[III, Exercise 8.3]{RH} we obtain:
\[
R^i\pi_*(\co_{\tilde{M}} \otimes \pi^*\cl) \cong R^i\pi_*\co_{\tilde{M}} \otimes \cl
\]
which implies that $R^i\pi_*\pi^* \cl = 0$ for any $i > 0$. 
The vanishing of the higher direct image of $\pi^* \cl$ and the isomorphism $\pi_*\pi^* \cl \cong \cl$ imply that $H^i(\tilde{M}, \pi^*\cl) \cong H^i(M,\cl)$.
\end{proof}

\begin{rem}
In the course of this paper, lemma~\ref{H^i inv sheaf} will be frequently used without referring to it.
For example, the proof of the vanishing of $\cl_X(d;m_1,\ldots,m_r)$ on $X$ will immediately give the vanishing of $\cl_Y(d;m_1,\ldots,m_r;\{ 0 \}_{i \in I})$ on $Y_I$.
\end{rem}

\section{The algorithm}\label{algorithm}

In this section we provide an algorithm for reducing any linear system
$\cl$ to one of four types of standard systems, for which we can determine its dimension explicitely. 
To obtain the speciality of $\cl$ it is then sufficient to
compare $\edim(\cl)$ and $\dim(\cl)$. 

\subsection{Reducing to a standard class}
Given a linear system $\cl$ we describe an algorithm for finding a new linear system $\cl'$ which is in standard form and such that $\dim \cl = \dim \cl'$.

\vspace{5mm}

{\bf Input:} $(d,m_1,\ldots,m_r)$.

\begin{enumerate}[]\itemsep 5mm \leftskip 5mm

\item {\bf Sort} the vector $(m_1,\ldots,m_r)$ in decreasing order.
\item {\bf While} ( $2d < m_1+m_2+m_3+m_4$\hspace{3mm} {\bf and}\hspace{3mm} $d > m_1$ )
\item \{

\vspace{5mm}
\begin{enumerate}[]\itemsep 5mm \leftskip 9mm
\item $k:= 2d-m_1-m_2-m_3-m_4.$
\item $d := d+k$\hspace{5mm} $m_i := \max(m_i+k,0)$\hspace{2mm}  ( $i=1\ldots 4$ )
\item {\bf Sort} the vector $(m_1,\ldots,m_r)$ in decreasing order.
\end{enumerate}
\item \}
\end{enumerate}
{\bf Output:} $(d,m_1,\ldots,m_r)$. 
\vspace{5mm}

Observe that if $d < m_1$ the system is empty. In this case the algorithm exit immediately from the main cycle and returns us a list which corresponds to an empty system.

The condition $2d < m_1+m_2+m_3+m_4$ has to be satisfied in order to decrease the degree of the system by means of a Cremona transformation. As mentioned in the previous paragraph, applying the Cremona transformation does not change the dimension of the linear system.

Finally, observe that if $m_i < 0$ then $-m_i E_i$ is in the base locus of $\cl$ so that
\[
	\dim \cl = \dim \cl + m_iE_i
\]	
and this justifies our redefinition of $m_i$.

After applying this algorithm, you either obtain that $\dim(\cl)= -1$ or $\dim \cl = \dim(\cl')$ where $\cl'$ is a standard class with $d \geq m_1\geq 0$.

\vspace{5mm}\subsection{Three types of standard classes}~

In what follows the symbol $\cl_X$ will always denote a non-empty standard linear system of the form $\cl_X(d;m_1,\ldots,m_r)$.
Associated to this system we will consider also its restriction to $Q_r$ which will be denoted by $\cl_Q$ and is of the form 
$\cl_Q(d,d;m_1,\ldots,m_r)$. Observe that also this system is standard. 
We will use the symbol $\cl$ to denote any one of $\cl_X$ or $\cl_Q$.
With this notation we have that:
\[
C\cdot\cl = 4d -m_1-\dots-m_r,
\]
where $C\in |-K_{Q_r}|$ is the strict transform of the smooth quartic curve through the $p_i$'s.

\subsubsection{{\bf Step 1: Splitting up according to $t$}}\label{step 1}~

\begin{claim*}
If $C\cdot\cl\geq 1$, then $C$ is not in the base locus of $\cl_Q$.\nvs
\end{claim*}

\begin{proof}
Since $\cl_Q$ is standard and $\cl_Q \cdot K_{Q_r} \geq 1 $, the result follows from Proposition~\ref{har} (because the blown up points are general on $C$).
\end{proof}

\begin{claim*}\label{cl3}
If $C\cdot\cl\leq 0$ and $\cl_X$ is not of the form $\cl_X(2m ; m^8 , m_9 , \ldots ,m_r)$, define
\begin{align*}
b &:= \max \{ i\ |\  4d - m_1 - \cdots - m_i + m_i (i-8) \geq 1\hspace{5mm} {\rm and }\hspace{5mm}  9 \leq i \leq r \}, \\
t &:= \left\lceil \frac{m_1 + \cdots + m_b -4d +1 }{b-8} \right\rceil
\end{align*}
Then $tC$ is contained in the base locus of $\cl$.\nvs
\end{claim*}

\begin{proof}
Note that $r$ can not be smaller than $8$ because in this case,
since $\cl_X$ is standard, we would have that either $c\geq 1$ or the system is of the form $\cl_X(2m;m^8)$. This implies that $b$ is well defined.

We proceed by observing that since $\cl_Q$ is standard and $\cl_Q \cdot K_{Q_r} \geq 0$ then, by Proposition~\ref{har}
we have that $C\subseteq\bs(\cl_Q)$. Removing $C$ from the system we obtain: 
\[
\cl_Q = C + \cl^1_Q,
\]
where $\cl^1_Q$ is still standard and $\cl^1_Q \cdot  K_{Q_r} = \cl_Q \cdot  K_{Q_r} - (r-8)$.
We repeat this procedure until the intersection of $C$ with the residual system $\cl^{t}_Q$ is positive or $t = m_r$.

If $b = r$, then, taking $t$ as in the statement of the claim,
$\cl^{t}_Q \cdot K_{Q_r} < 0$ and $\cl^{t-1} \cdot K_{Q_r} \geq 0$,
which implies that $tC\subseteq \bs(\cl_Q)$.

If $b < r$, then $m_r C$ is contained in the base locus of $\cl_Q$. Consider now the decomposition $\cl_Q = m_r C + \cl^{m_r}_Q$ and 
let $r':= \max\{ i\ |\ m_i - m_r > 0 \}$, then the class of
$\cl_Q^{m_r}$ is the pull-back of a class on $Q_{r'}$  so that we can repeat the above arguments.

In case $b=r'$, proceeding as before, we obtain that, for
\begin{align*}
t & = m_r + \left\lceil\frac{(m_1-m_r) + \cdots + (m_{r'}-m_r) - 4d + 8 m_r + 1 }{r'-8}\right\rceil \\
  & = \left\lceil\frac{m_1 + \cdots + m_{r'} - 4d + 1 }{r'-8}\right\rceil,
\end{align*}
the divisor $tC$ is contained in the base locus of $\cl_Q$.

Since $b \geq 9$ exists, after applying these arguments a sufficient number of times, our procedure comes to an end and we obtain that $tC$ is contained in the base locus of $\cl_Q$ for the claimed value of $t$.
\end{proof}

\begin{concl}
We distinguish the following cases for $\cl_X$:
\begin{enumerate}[{\rm 1.}]\leftskip 7mm\itemsep 3mm
\item $C\cdot\cl_X\geq 1$, i.e. $t=0$.
\item $\cl_X(2m ; m^8 , m_9 , \ldots ,m_r)$, i.e. $t = m$.
\item $C\cdot\cl_X\leq 0$ and $\cl_X \neq  \cl_X(2m ; m^8 , m_9 , \ldots ,m_r)$,
	i.e. $t$ is as in claim~\ref{cl3}.
\end{enumerate}
\end{concl}

Note that we haven't proved that $t=m$ in case (2), but the above arguments will still work
to obtain that $m_9 C$ is contained in the base locus of $\cl_Q$. Since the residue
class is then $\cl_Q( 2m',2m'; m'^8)$, with $m'= m - m_9$,
the fact that $t = m$ follows from~\cite[Proposition 1.2]{bh1}.

\subsubsection{{\bf Step 2 : Reducing to the case $t\leq m_r$}}\label{step 2}~

In case (1) there is nothing to be done since $t=0$.

In case (2), since the system is standard, we have that $m \geq m_i$
for any $i$. The fact that $m C$ belongs to the base locus of $\cl_X$ implies that 
\[
\dim \cl_X(2m ; m^8 , m_9 , \ldots ,m_r) = \dim \cl_X(2m ; m^r),
\]
so it is sufficient to determine the dimension of the last system.

In case (3) we have that $t \leq m_r$ if and only if $b=r$. 
On the other hand, if $b < r$ then $m_{b+1} < t \leq m_b$ and, since $tC$ is contained in the base locus of $\cl_X$ we have:
\[
\dim \cl_X(d;m_1,\ldots,m_r) = \dim \cl_X(d;m_1,\ldots,m_b,t^{r-b}).
\]
As before we need just to determine the dimension of the last linear system.

\subsubsection{{\bf Step 3: Reducing to the case $d \geq m_1 + t$}}\label{step 3}~

The previous part allows us to limit our study to the case:
\[
\cl_X \neq  \cl_X(2m ; m^8 , m_9 , \ldots ,m_r),\hspace{1cm}
C\cdot\cl\leq 0\hspace{1cm} {\rm and}\hspace{1cm} t \leq m_r.
\]

In what follows we will adopt the notation: 
\[
t_1 = \max\{0,m_1+m_2-d\}\hspace{1cm} {\rm and}\hspace{1cm} t_i := \max\{0,m_1+m_i-d\}\hspace{5mm} {\rm for} \hspace{5mm}i=2,\ldots,r.
\]

\begin{claim}\label{cone in bs}
Assume that $d < m_1 +t$, then $r \geq 10$, the divisor 
$M \in \cl_X(3;3,1^{r-1})$ is a cone contained in the base locus of $\cl_X$ and $t_r > 0$.

\end{claim}

\begin{proof}
Since we reduced to consider the case $t\leq m_r$, the positivity of $t_r$ follows immediately from the hypothesis.
Observe that $t$ must be positive and that $\cl_X$ can not be of the form $\cl_X(2m ; m^8 , m_9 , \ldots ,m_r)$. As noticed before this implies that $r\geq 9$. Observe that the equality does not hold, since  we in this case $t = m_1 + \cdots + m_9 -4d +1$ and 
\[
m_1 + t - d = (m_1 + \cdots + m_4-2d)+(m_1 +m_5 \cdots + m_7-2d)+(m_8+m_9-d) +1,
\]
where the first two terms in parentheses are non-positive and the third is negative because $\cl_X$ is in standard form. This implies that $r$ is at least $10$.

The linear system of $\p^3$ given by $\cl_3(3;3,1^9)$ contains a unique divisor which is a cone over a plane cubic. Since 
\[
C\cdot\cl_X(3;3,1^9) = 0
\]
we have that $C$ is contained in the base locus of this system by Proposition~\ref{har}. This implies that 
$\dim \cl_X(3;3,1^{r-1}) = 0$ if $r \geq 10$ and the unique element $M$ is the strict transform of the cone with vertex $p_1$ and base curve $B$. 

Let $D\in\cl_3(d;m_1,\ldots,m_r)$ be a general element of the system that we are considering. Observe that any line $l$ through $p_1$ and $q\in B$ has an intersection multiplicity with $D$ at least
\[
m_1+t > d.
\]
This means that the strict transform of $l$ is contained in the base locus of $cl_X$ and this implies that $M$ is contained in the base locus of $\cl_X$.
\end{proof}

\begin{claim}\label{L - cone}
Assume that $d < m_1 +t$, then the system $\cl'_X := \cl_X - \cl_X(3;3,1^{r-1})$ is almost standard or empty.
\end{claim}

\begin{proof}
In order to see that $\cl'_X$ is almost standard, it is sufficient to check that the multiplicities $m_1 -3, m_2-1,\ldots,m_r-1$  are non negative and that $2d - 6$ is bigger or equal to the sum of the biggest four multiplicities.
Observe that if the inequality 
\[
m_1 -3 \geq m_5 -1
\]
holds, then $m_1-3$ belongs to the set of biggest four multiplicities 
and this gives the thesis.
Assume that $m_1 - 3 < m_5 -1$ then we would have $m_1 -1 \leq m_5 \leq m_1$. The fact that $\cl_X$ is standard together with claim~\ref{cone in bs} imply that
\[
2d \geq m_1 + \cdots + m_4\hspace{1cm} {\rm and}\hspace{1cm} d < m_1+m_5.
\]
By substituting the two possible values for $m_5$ we obtain a contradiction.
\end{proof}

Summarizing the two claims we see that if $d < m_1 + t$ then, either $cl_X$ is empty or it is one of the following:
\begin{enumerate}[{\rm i.}]\leftskip 7mm\itemsep 3mm
\item $\cl_X(3;3,1^r)$ with $r \geq 10$ and $\dim \cl_X = 0$;
\item $M + \cl'_X$ where $\cl_X'$ is almost standard and $\dim\cl_X = \dim\cl'_X$.
\end{enumerate}

In case (i) we know the dimension of $\cl_X$.
For case (ii), after reordering the multiplicities of $\cl_X'$ we obtain a new system $\cl''_X$ which is standard.
To determine the dimension of $\cl''_X$ we repeat the procedure starting from Step 1 of section \ref{step 1}.

{\bf Conclusion}~\label{concl 3 types}

Using the above Steps, we see that we are reduced to determing the dimension of $\cl_X$ in the following cases:
\begin{enumerate}[{\rm 1.}]\leftskip 7mm\itemsep 3mm
\item $C\cdot \cl_X\geq 1$, i.e. $t=0$.
\item $\cl_X =  \cl_X(2m ; m^r)$ with $r \geq 8$, i.e. $t = m$.
\item $1+m_r(8-r)\leq C\cdot\cl_X\leq 0$ and $d \geq m_1 + t$ where $t = \displaystyle\left\lceil\frac{C\cdot\cl +1 }{r-8}\right\rceil \leq m_r$ and $\cl_X\neq\cl_X(2m ; m^8 , m_9 , \ldots ,m_r)$,

\end{enumerate}

\section{Determining the dimension of $\cl$ for case (1)}

We begin by defining $I$ to be set of indices $i$ such that $t_i>0$.
Observe that if $I$ is non-empty then, since $\cl_X$ is standard, it must be one of these two types:
\[
\{2,3,\ldots,s-1,s\}\hspace{1cm} {\rm or}\hspace{1cm} \{1,2,3\}.
\]

The following theorem shows that $h^1(X,\cl_X)$ is a function 
of the numbers
\[
t_i = - l_i\cdot\cl_X
\]
where the $l_i$'s are the strict transforms of lines indexed by $I$.
This fact has already been proved in~\cite{DL8} for the case $r \leq 8$. We will however make the proofs in this paper self-contained.

\begin{theor}\label{t=0}
Let $\cl_X$ be a non-empty standard system with $C\cdot\cl_X\geq 1$,
then 
\[
\dim \cl_X = v(\cl_X) + \sum_{i=1}^r \binom{t_i +1}{3}.
\]
\end{theor}

\begin{proof}
Without loss of generality, we may assume that $m_r >0$.
This means our linear system $\cl_X$ satisfies the following conditions:
\begin{enumerate}[{\rm 1.}]\leftskip 7mm\itemsep 3mm
	\item $2d \geq m_1 + m_2 + m_3 + m_4$
	\item $d \geq m_1 \geq \ldots \geq m_r > 0$
	\item $t_1 := \max \{ 0 , m_2 + m_3 -d \}$ 
		and $t_i := \max \{ 0, m_1 + m_i - d \}$ for $i = 2,\ldots,r$
	\item $4d \geq m_1 + \cdots + m_r +1$
\end{enumerate}

Consider the blow-up $\pi_I: Y_I\map X$ along the $l_i$, with $i \in I$ and let $\cl_Y$ be the complete linear system
$\cl_Y(d;m_1,\ldots,m_r;\{t_i\}_{i\in I})$ (and its corresponding invertible sheaf) on $Y_I$.
Since, for all $i \in I$, the curve $t_il_i$ belongs to the base locus of $\cl_X$, we have that $\dim\cl_X = \dim\cl_Y$.

Now, using lemma's~\ref{c1,c2} and \ref{int Y_I}, an easy but tedious calculation shows
that
$$
\cx(Y,\cl_Y) = \cx(X,\cl_X) + \sum_{i=1}^r \binom{t_i +1}{3}.
$$
So, in order to prove theorem~\ref{t=0}, it is sufficient to show that for $i\geq 1$
\begin{equation}\label{eq h^i}
h^i(Y_I,\cl_Y) =  0.
\end{equation}

Denote by $\cs_r := \cl_Y(2; 1^r)$ and consider the exact sequence associated to a smooth $\tilde{Q} \in \cs_r$
$$
\xymatrix@1@C=40pt{
0 \ar[r] & \cl_Y - \cs_r \ar[r] & \cl_Y \ar[r] & \cl_Y \otimes \co_{\tilde{Q}} \ar[r] & 0.
}
$$

Now, since $\tilde{Q}$ does not intersect the exceptional divisor $F_i$ of $l_i$,
we see that $\cl_Y \otimes \co_{Q_r} = \cl_Q(d,d;m_1,\ldots,m_r)$
is standard. It then follows from condition (4) and \cite[Theorem 1.1]{bh1} that for $i\geq 1$
\[
h^i(\cl_Y \otimes \co_{\tilde{Q}}) = 0.
\]

In order to obtain~(\ref{eq h^i}), we obviously need the vanishing of the higher cohomology groups of
$\cl_Y - \cs_r = \cl_Y(d-2;m_1-1,\ldots,m_r-1,\{t_i\}_{i \in I})$.
We begin by distinguishing two cases according to the values of $r$.
\begin{claim}\label{r<=3}
If $r\leq 3$ then $h^i(Y,\cl_Y) = 0$ for all $i\geq 1$.
\end{claim}

Observe that as long as $r\geq 4$, the system $\cl_Y - \cs_r$ satisfies conditions (1)-(4) unless one of the following occurs:
\begin{enumerate}[{\rm (a)}]\leftskip 7mm\itemsep 3mm
	\item $d > m_1$ and $m_r = 1$
	\item $d = m_1$ and $m_r > 1$
	\item $d=m_1$ and $m_r =1$
\end{enumerate}

Recall that $s$ is defined to be the maximum of $I$ and observe that if $d = m_1$ then $s=r$. 

So, if $s < r$, the system $\cl_Y - \cs_r$ satisfies conditions (1)-(4), unless $m_r =1$.
However, in this case, we can consider it as
a linear system on $Y_{r'I}$, where $r' := \max\{ i\ |\ m'_i >0\}$.
If we can prove the vanishing of the cohomology groups of $\cl_Y - \cs_r$ on
$Y_{r'I}$, then, because of lemma~\ref{H^i inv sheaf}, this implies the vanishing of 
the cohomology groups on $Y_I$.
Note that, $r'\geq s$.

Moreover, if $r=s$, then case (a) cannot occur, because this would imply that
$m_1 + m_r - d = m_1 + 1 - d \leq 0$, which contradicts $t_r > 0$.

In any case, if $d' > m'_1$, then we can consider an exact sequence like before
(using $\cs_{r'}$ in stead of $\cs_r$ if $m_r = 1$),
and, using arguments as before,
we can reduce to proving the vanishing of the higher cohomology groups
of $\cl_Y - 2 \cs_r$ (or $\cl_Y - \cs_r - \cs_{r'}$ if $m_r=1$).
Obviously, this procedure can be repeated
untill we are reduced to proving the vanishing of the higher cohomology groups for a class 
$\cl_Y(d; m_1, \ldots, m_r; \{t_i\}_{i \in I} )$ satisfying conditions (1)-(4), and with either $r \leq 3$ or $r = s \geq 4$ and $d = m_1$.
Note that $t_i$ and $I$ occurring in $\cl_Y$ are in fact the same as the original ones.

We thus end the proof with the following:

\begin{claim}\label{d=m1}
Assume that $\cl_Y = \cl_Y(m_1; m_1 ,\ldots,m_r; \{t_i\}_{i \in I})$ satisfies conditions (1)-(4) with $s=r \geq 4$, then $h^i(Y,\cl_Y-\cs_r) = 0$ for all $i\geq 0$.
\end{claim}

\end{proof}

\begin{proof}[Proof of Claim~\ref{r<=3}]
This can be regarded as the ``toric case'' by putting the three points $p_1, p_2, p_3$ in $(1:0:0:0), (0:1:0:0), (0:0:1:0)$ and observing that $X$ and $Y_I$ are toric varieties. The evaluation of the dimension of $\cl_X$ can be worked out by counting the monomials with multiplicities $m_i$ at these points. In the same way the cohomology groups of $\cl_X$ and $\cl_Y$ can be found by purely combinatorial methods and in particular we have that $h^i(Y,\cl_Y) = 0$ for all $i\geq 1$.

\end{proof}

\begin{proof}[Proof of Claim~\ref{d=m1}]
By hypothesis we have that $t_i = m_i$ for $i =2,\ldots,r$.
Now let $a_i := m_i - 1$ for $i=1,\ldots,r$ so that
$$\cl_Y - \cs_r = \cl_Y(a_1-1; a_1 ,\ldots, a_r ; a_2 + 1, \ldots , a_r +1)$$
with $a_1 \geq \cdots \geq a_r \geq 0$ and denote this class by $\ca_r$.

We will now prove by induction on $b$ that for any $i\geq 0$
$$
h^i(\ca_b) = 0,
$$
for all $b \geq 1$ and $a_1 \geq \cdots \geq a_b \geq 0$.

An easy calculation, using lemma's~\ref{c1,c2} and \ref{int Y_I},
shows that $\cx(\ca_b)=0$.
Moreover, $h^0(\ca_b) = 0$ because the degree $a_1-1$ is less than the biggest multiplicity $a_1$. In this way it is enough to prove the vanishing of the first and second cohomology group.

If $a_1 = 0$, then $\ca_1$ is the system $| - H |$ so that $h^1(\ca_1) = h^2(\ca_1) = 0$ (see e.g. \cite[III, Theorem 5.1 (b), p. 225]{RH}).

Now assume $b=1$ and $a_1 > 0$, so that $\ca_1 = \cl_X(a_1-1; a_1)$.
Note that, because $b=1$, we do not blow up lines, so we work on $X_1$.
Let $W \in \cw = \cl_X(1;1)$ and consider the exact sequence
$$
\xymatrix@1@C=40pt{
0 \ar[r] & \ca_1 - \cw \ar[r] & \ca_1 \ar[r] & \ca_1 \otimes \co_{W} \ar[r] & 0.
}
$$
Since $\ca_1 - \cw = \cl_X(a_1-2; a_1-1)$, our induction hypothesis
implies that $h^i(\ca_1 - \cw) =0$ for $i=0,1,2$.
On the other hand, because $\ca_1 \otimes \co_{W} = \cl_2(a_1-1;a_1)$, one easily checks that $h^i(\ca_1 \otimes \co_{W}) = 0$ for all $i\geq 0$.

Next, assume that $b > 1$ and that the statement is true for $b' \leq b-1$.
On $Y_{bI}$, with $I=\{2,\ldots,b\}$, consider the exact sequence
$$
\xymatrix@1@C=40pt{
0 \ar[r] & \ca_{b-1} - F_b \ar[r] & \ca_{b-1} \ar[r] & \ca_{b-1} \otimes \co_{F_b} \ar[r] & 0.
}
$$
Recall that, since $F_b\cong\p^1\times\p^1$, the restriction of $\ca_{b-1}$ to $F_b$ is given by:
\[
\ca_{b-1} \otimes \co_{F_b} \equiv \co(-1,0).
\]
A standard argument shows that $h^i(F_b,\co(-1,0)) = 0$ for all $i\geq 0$ and this, together with the vanishing of $h^i(\ca_{b-1}) = 0$ ( by induction ) implies that $h^i(\ca_{b-1} - F_b)=0$ for $i\geq 0$.
But 
\[
\ca_{b-1} - F_b = \ca_b\hspace{1cm} {\rm with} \hspace{1cm} a_b = 0.
\]

Again we will use induction, now on $a_b$,
so we assume that the statement is true for all $a'_b \leq a_b -1$.
Let $\ca'_b = \cl_Y(a_1-1; a_1 ,\ldots, a_{b-1}, a_b-1 ; a_2 + 1, \ldots , a_{b-1} +1, a_b)$
and consider the exact sequences
$$
\xymatrix@1@C=40pt{
0 \ar[r] & \ca'_{b} - E_b \ar[r] & \ca'_{b} \ar[r] & \ca'_{b} \otimes \co_{E_b} \ar[r] & 0
}
$$
and
$$
\xymatrix@1@C=30pt{
0 \ar[r] & \ca'_{b} - E_b - F_b \ar[r] & \ca'_{b} - E_b \ar[r] & 
		(\ca'_{b} - E_b) \otimes \co_{F_b} \ar[r] & 0.
}
$$
Since $\ca'_{b} \otimes \co_{E_b} = \cl_2(a_b-1;a_b)$, we have that $h^i(\ca'_{b} \otimes \co_{E_b})=0$ for $i\geq 0$.
So, because of our induction hypothesis, the cohomology of the first exact sequence implies
$h^i(\ca'_{b} - \co_{E_b})=0$ for $i=0,\ldots,3$.
On the other hand $(\ca'_{b} - E_b) \otimes \co_{F_b} \equiv \co(-1,a_b)$, and one easily checks that
$h^i((\ca'_{b} - E_b) \otimes \co_{F_b})=0$ for $i\geq 0$.
So, from the cohomology of the second exact sequence and the fact that $\ca'_{b} - E_b - F_b = \ca_b$ we conclude that
$h^i(\ca_{b})=0$ for $i\geq 0$.
\end{proof}

\section{Determining the dimension of $\cl$ for case (2)}

As seen before, the dimension of $\cl_X(2m;m^8,m_9,\ldots,m_r)$ is equal to the dimension of $\cl_X(2m;m^r)$ and this number is evaluated in the following.

\begin{theor}\label{2m,m^r}
Let $r \geq 8$  then
$\dim \cl_3 (2m ; m^r) = m$.
\end{theor}

\begin{rem}
The virtual dimension of the system $\cl_3 (2m ; m^8)$ is equal to $m$ only if $r=8$. This means that for bigger values of $r$ the system is special.
\end{rem}

\begin{proof}
The statement is trivial for $m=0$, so assume $m$ to be positive.
Because of \cite[Theorem 6.2 (2)]{DLB}, we know that
the base locus of the linear system $\cl_3 (2m ; m^8)$ on $X_8$,
is $mC$.
This implies that $\dim \cl_3 (2m ; m^8) = \dim \cl_3 (2m ; m^r)$ for all $r \geq 8$
(because all the $p_i$'s lie on $C$).
And according to \cite[Theorem 5.1]{DL8}, we then obtain the statement.
\end{proof}

\begin{rem}
Let us just note that the techniques used throughout this paper can yield an alternative prove of
theorem~\ref{2m,m^r}. 
More precisely, we can consider the blowing up $\overline{Y}_{\! r}$ of $X$ along $C$,
and the linear system $\tilde{\cl}_Y := \tilde{\cl}_Y(2m;m^r;m)$ on $\overline{Y}_{\! r}$.
Using induction on $m$ and cohomology of the exact sequence
$$
\xymatrix@1@C=40pt{
0 \ar[r] & \tilde{\cl}_Y - \cs_r \ar[r] &  \tilde{\cl}_Y \ar[r] & 
		\tilde{\cl}_Y \otimes \co_{\tilde{Q}} \ar[r] & 0,
}
$$
with $\tilde{Q} \in \cs_r := \tilde{\cl}_Y(2;1^r;1)$,
one can obtain that
$\dim (\tilde{\cl}_Y) = m$.
And, since obviously $\dim \cl_3 (2m ; m^r) = \dim \tilde{\cl}_Y$,
this proves the theorem.
\end{rem}

\section{Determining the dimension of $\cl$ for case (3)}

We recall that $t_i := \max \{ 0, m_1 + m_i - d \}$, for $i=2,\ldots
r$ is equal to the opposite of the intersection of the strict transform of $l_i$ with $\cl_X$ ( if this intersection is negative ).
The same is true for $t_1 := \max \{ 0 , m_2 + m_3 -d \}$, where the line is $l_1$ through $p_2$ and $p_3$.
Finally, the number $t :=  \left\lceil\frac{C\cdot\cl+1 }{r-8}\right\rceil $ can be described as  
\[
t := \max\{i\in\mathbb{N}\ |\ C\cdot (\cl_X-iC) \leq 0\}.
\]

The next theorem shows that the speciality of $\cl$ only comes from lines $l_i$ for which $t_i >0$ and from the curve $C$.

\begin{theor}\label{t>0,a}
Let $\cl(d;m_1,\ldots,m_r)$ be a standard class which is not of type (2). Assume that
\[
1+m_r (8-r)\leq C\cdot\cl\leq 0
\]
and $d \geq m_1 + t$, then 
\[
\dim \cl = v(\cl) + \sum_{i=1}^r \binom{t_i +1}{3}+(r-8)\binom{t+1}{3} + n \binom{t+1}{2},
\]
where $n\leq r-9$ is a non-negative integer such that 
$n = (8-r)(t-1)+C\cdot\cl$.
\end{theor}

\begin{proof}
Note that, since $\cl_X$ is standard and $C\cdot\cl\leq 0$ we have that $r$ is at least $9$. From the assumptions, it follows that $0 < t \leq m_r$.
Define $I$ to be the set of $i$ such that $t_i$ is positive and consider the system $\tilde{\cl}_Y$ defined on $\tilde{Y}_I$ of the form:
\[
\cl_Y(d ; m_1 , \ldots , m_r; \{ t_i\}_{i \in I}, t).
\]
Since $\sum_{i \in I} t_i l_i + t C$ belongs to the base locus of $\cl_X$,
we have that $\dim \cl = \dim \tilde{\cl}_Y$.
The linear system $\tilde{\cl}_Y$ thus satisfies the following conditions
\begin{enumerate}[{\rm 1.}]\leftskip 7mm\itemsep 3mm
	\item $2d \geq m_1 + m_2 + m_3 + m_4$
	\item $d \geq m_1 \geq \ldots \geq m_r > 0$
	\item $t_1 := \max \{ 0 , m_2 + m_3 -d \}$ 
		and $t_i := \max \{ 0, m_1 + m_i - d \}$ for $i = 2,\ldots,r$
	\item $t = \left\lceil\dfrac{C\cdot\cl +1 }{r-8}\right\rceil 
			= \dfrac{C\cdot\cl + r-8 - n }{r-8}$, $0 < t \leq m_r$
	\item $d \geq m_1 + t$.
\end{enumerate}
On the other hand, using lemma's~\ref{c1,c2} and \ref{int Y_I}, 
an easy but tedious calculation shows that
\[
\cx(Y,\cl_Y) = \cx(\tilde{Y},\tilde{\cl}_Y) + (r-8)\binom{t+1}{3} + y \binom{t+1}{2}.
\]
So, in order to prove theorem~\ref{t>0,a}, it is sufficient to show that for any $i\geq 1$
\begin{equation}\label{h^i 3a}
h^i(\tilde{Y}_I,\tilde{\cl}_Y) = 0.
\end{equation}

Denote by $\cs_r := \cl_Y(2; 1^r ; \{0\}_{i \in I} ; 1)$ and consider a general element $\tilde{Q} \in \cs_r$ from which we have the exact sequence
$$
\xymatrix@1@C=40pt{
0 \ar[r] & \tilde{\cl}_Y - \cs_r \ar[r] & \tilde{\cl}_Y \ar[r] 
	& \tilde{\cl}_Y \otimes \co_{\tilde{Q}} \ar[r] & 0.
}
$$

By abuse of notation, let $C$ also denote the anticanonical curve on 
$\tilde{Q} \subset \tilde{Y}_I$,
then $F|_{\tilde{Q}} = C$.
So, since $\tilde{Q}$ and $F_i$ are disjoint,
we have that the restriction of $\tilde{\cl}_Y$ to $\tilde{Q}$ is:
\[
\tilde{\cl}_Y \otimes \co_{\tilde{Q}}  = \cl_Q(d-2t,d-2t;m_1-t,\ldots,m_r-t).
\]
Observe that this system is standard because $d \geq m_1 +t$ and $m_r \geq t$.
On the other hand, condition~(4) implies that
$\cl_Q \cdot K_{\tilde{Q}} \leq -1$.
So it follows from \cite[Theorem 1.1]{bh1} that 
$h^i(\tilde{\cl}_Y \otimes \co_{\tilde{Q}}) = 0$ for any $i\geq 1$.

In order to obtain~\eqref{h^i 3a}, we obviously need the vanishing of the higher cohomology
groups of $\tilde{\cl}_Y - \cs_r$, which is of the form:
\[
\tilde{\cl}_Y(d-2;m_1-1,\ldots,m_r-1,\{t_i\}_{i \in I};t-1).
\]
First, let us check when 
$\cl_Y (d' ; m'_1 , \ldots , m'_r, \{t'_i\}_{i \in I};t') = \tilde{\cl}_Y - \cs_r$ 
satisfies conditions (1)-(5).
Observe that conditions (1), (3) and (5) are always satisfied.

(2) is satisfied unless $m_r =1$ and in this case $t = 1$.

(4) Since 
    \begin{align*}
	t' & = \frac{m_1 + \cdots + m_r -4d + r-8 - n }{r-8} -1 \\
		& = \frac{m'_1 + \cdots + m'_r -4d' + r-8 - n }{r-8},
    \end{align*}
	we obtain $0 < t' \leq m'_r$ unless $t=1$.

In this way we see that $t >1$, then the system satisfy all the conditions. This means that, if $t' \geq 1$, we can consider an exact sequence as above,
and, using arguments as before, we can reduce to proving the vanishing of the higher cohomology groups of
$\tilde{\cl}_Y - 2\cs_r$.
Obviously, this procedure can be repeated $t$ times,
or thus untill we are left with proving the vanishing of the higher cohomology groups of 
$$
\tilde{\cl}_Y - t\cs_r = \tilde{\cl}_Y(d-2t; m_1-t,\ldots,m_r-t; \{t_i\}_{i \in I} ; 0 ).
$$
Since $t=0$, we only need to prove the vanishing of the higher cohomology groups on $Y_I$ of the
class
$\tilde{\cl}_Y(d-2t; m_1-t,\ldots,m_r-t; \{t_i\}_{i \in I} ) $,
which we will denote as 
$\cl_Y$.
Now, since 
\[
C\cdot\cl_Y -1 = 4d - 8t - m_1 - \cdots - m_r + rt - 1 
		 = r - 9 - n \geq 0,
\]
we have that the class $\cl_Y$ satisfies the conditions of theorem~\ref{t=0} and we conclude that $h^i(\cl_Y) =0$ for $i\geq 1$.
As mentioned before, this is enough to prove the theorem.
\end{proof}

\bibliographystyle{alpha}

\end{document}